\documentclass{iet-ell}

\usepackage{algorithmic}
\usepackage{algorithm}
\usepackage{cuted} 
\usepackage{tikz}
\usepackage{graphicx}
\usepackage{xcolor} 
\usetikzlibrary{decorations.pathmorphing, positioning, shapes.geometric}
\usetikzlibrary{decorations.markings}
\usetikzlibrary{angles, quotes}
\setcopyright{open}
\ietVolume{00}
\ietYear{2021}
\ietDoi{ell2.10001}

\received{10 January 2021}
\accepted{4 March 2021}

\begin{document}

\title{Multistatic anisotropic travel-time imaging as a tensor tomography problem}

\author[af0]{Naeem Desai}
\orcid{0000-0002-2071-1788}
\author[af1]{Oliver Graham}
\orcid{0009-0006-9451-1152}
\author[af1]{Sean Holman}
\orcid{0000-0001-8050-2585}
\author[af1]{William R.B. Lionheart}
\orcid{0000-0003-0971-4678}

\orcid{1234-5678-9012-3456a}
\affil[af0]{DSTL, UK}
\affil[af1]{Department of Mathematics, The University of Manchester, UK. Note author order is alphabetic}

\corresp{Email: bill.lionheart@manchester.ac.uk}

\begin{abstract}%
Travel-time imaging problems seek to reconstruct an image of reflectivity of a scene by measuring travel time (and amplitude, phase) of electromagnetic or acoustic signals, such as radar and sonar.  Multistatic, in this context, means that the transmitters and receivers need not be co-located. The reflectivity is anisotropic if it depends on direction, and in the multistatic case this means incoming and outgoing direction. Travel-time problems can be formulated as generalized Radon transforms of integrals over isochrones, in the planar case ellipses with transmitter and receivers at foci.  In a simplified case where transmitters and receivers are distant from the scene, isochrones can be approximated by straight lines. We relate this to tensor ray transforms, specifically the longitudinal ray transform of Sharafutdinov, and discuss the implication of its known null-space. In the volumetric case isochrones are spheroids and we relate the problem to the normal Radon transform of tensor fields.
\end{abstract}

\maketitle%

\section{Introduction}\label{sec:Introduction}

Many inverse problems in imaging can be considered multistatic in that for each transmitter (source) position signals are detected by  multiple receivers (detectors).  In X-ray tomography the source position moves and a pixel array of detectors captures a projection image. In electrical impedance tomography the sources and detectors are simply electrodes. In magnetic induction imaging methods, including imaging metal detectors, one coil transmits while the electromagnetic field is detected by multiple coils \cite{ma2006imaging}. In passive  radar transmitters of opportunity, on vehicles for example, are used that are not under the operators' control while multiple, possibly mobile, receiver stations are used to detect reflected signals. Active multistatic radar is in its infancy, but rapid improvements in position and time measurements, including quantum accelerometers and chip-scale electronic clocks, make it conceivable at least to have a large collection of mobile transmitters and receivers with known positions and synchronized clocks. Similar considerations apply to acoustic imaging such as sonar. 

In the most detailed model the governing equation would simply be the electromagnetic or acoustic wave equation, with spatially varying material properties.  The inverse problem for recovering these properties is non-linear and non-local, and would involve iterative solution of the forward problem with a numerical solver (such as finite  element or finite difference). The forward model would incorporate all the relevant physics, including diffraction and multiple scattering.  This approach is termed full wave inversion and there are examples in medical ultrasound \cite{ali20242}, ground penetrating radar \cite{watson2016better}
  and seismic imaging. 
  Such methods are computationally intensive. 

In this paper we consider a simplified model with a lower computation cost.  In some cases it allows us to use known results in tensor tomography while it also suggests some future research that extends these ideas.  Essentially we use a two scale model.   We consider clusters of objects that are below the spatial resolution we expect that are treated as having an aggregate response to an incoming plane wave in direction $\xi_{\mathrm{in}}$ as a function of the outgoing
direction $\xi_{\mathrm{out}}$. The possibly complex reflectivity
  $f(x,\xi_{\mathrm{in}},\xi_{\mathrm{out}})$ 
 is a function of the unit vectors associated with these directions and the position $x$ on the larger length scale. In practice for a radar application we might consider a small tree, vehicle or house, in which there are multiple reflections on the small scale, with the overall multistatic radar scattering cross section summarized by $f$. Reversibility implies that
  $f(x,\xi_{\mathrm{in}},\xi_{\mathrm{out}}) = f(x,\xi_{\mathrm{out}},\xi_{\mathrm{in}})$. 
   
  We expect that as a function of $x$, $f$ may well be discontinuous, for example at the boundary of a type of vegetation.  A simple trihedral corner reflector would have the  property that within a range of incoming directions there is only a non-zero signal for outgoing directions close to the incoming direction. In that case $f$ would be an approximate delta function as a function of the angle between $\xi_{\mathrm{in}}$ and $\xi_{\mathrm{out}}$.  As the complexity of the object increases we would expect the dependence on directions to be generally smoother, but with some peaks. There is a marked absence of published multistatic radar scattering cross section data, but we hope that this paper will stimulate experimental work. 
  
  Let us consider the scene, or region of interest, to be defined by the ball $|x|< R_{\mathrm{ scene}}$, so the support of $f$ as a function of $x$ is contained in the scene.  We will be interested in the case where $x_T$ and $x_R$ lie outside the scene, and typically with $|x_T|, |x_R| \gg R_{\text{\tiny scene}}$. 

\section{Isochrones and generalized Radon transforms}


In 2 or 3 dimensions, let $x_T$, $x_R$ be the transmitter and receiver positions. A reflection from a point $x$ in space has a travel time $t$ when
\begin{equation}
	|x - x_T| + |x-x_R|=ct 
\end{equation}
where $c$ is the wave speed. We call this set of points $x$ an {\em isochrone} and it is  an ellipse with foci at $x_T$ and $x_R$ in the 2D case and similarly in 3D a spheroid with these foci, with a rotational axis of symmetry on the line joining transmitter and receiver. We assume that the received signal at $x_R$ will be the sum of the reflections determined by $f$ for all the points on the isochrone. Of course, this ignores any multiple reflections. We are also ignoring the antenna radiation pattern, and while this is appropriate in ``searchlight'' mode SAR it could easily be incorporated in the model. We must be careful of the length or area measure used for this average. Moon \cite{moon2014determination} considers the isotropic case where $f$ depends only on $x$, and takes the view that the time is measured with a fixed precision, taking the integral
\begin{equation}
	E(x_T,x_R,t)[f] = \int\limits_{|x - x_T| + |x-x_R|<ct } f(x) \mathrm{d} x 
\end{equation}
and defining the ellipse or spheroid generalized Radon transform as
\begin{equation}\label{eq:EtoR}
R(x_T,x_R,t)[f] = \partial E(x_T,x_R,t)[f] /\partial t
\end{equation}
which differs from the standard  integral over a curve or surface. Palamodov \cite{palamodov2011uniform} gives a uniform formula for inversion of generalized Radon transforms, provided they are given in a similar form of level sets with respect to a parameter. 

For the anisotropic case we extend this idea by including the directional dependence and define 
\begin{equation}
	E(x_T,x_R,t)[f] = \int\limits_{|x - x_T| + |x-x_R|<ct } f(x, \widehat{x-x_T}, \widehat{x-x_R}) \mathrm{d} x 
\end{equation}
where $\widehat{v} = v/|v|$ for any vector $v$,  the anisotropic ellipse or spheroid generalized Radon transform  $R$ as in  (\ref{eq:EtoR}).  As yet there is no  theory for the anisotropic case unless we take the approximation in the next section.

\section{Flat isochrone approximations}

The maximum and minimum (Gaussian) curvature of an ellipse (spheroid) occur on the minor and major axes.  Let $d=|x_R -x_T|/2$ be half the focal distance, and assuming $ct > 2d$  then we see the  maximum curvature in  the 2D case  is 
\begin{equation}
	\kappa_{\mathrm{max}} = \frac{4 \sqrt{ (ct)^2 /4 -d^2 } }{(ct)^2} 
\end{equation}
and the maximum Gaussian curvature in the 3D case is
\begin{equation}
	K_{\mathrm{max}} = \frac{16 \left(  (ct)^2 /4 -d^2\right) }{(ct)^4}.
\end{equation}
To neglect this curvature we need the scene under consideration to be much smaller than the maximum  radius of curvature $1/\kappa_{\mathrm{max}}$  ($1/K_{\mathrm{max}}$)
then  we can approximate  a small part of the isochrone by its tangent line (plane).  The normal  vector is in the direction
\begin{equation}
n=\frac{ \widehat{x-x_T}  + \widehat{x-x_R}}{2}, 
\end{equation}
which we can think of as the {\em average azimuthal direction} as it bisects the lines from $x$ to the transmitter and receiver.  Let $s = |x - (x_T +x_R)/2|$ be the distance between $x$ and the centre of the isochrone, then the tangent line (plane) is $y: (y - (x_T +x_R)/2) \cdot \hat{n} =s$. 

\section{Tensor Ray transforms}

The normal Radon Transform of a symmetric rank $k$ tensor  field $F$ with components $f_{i_1,....i_k}$  then for each line (plane) $x \cdot \hat{n} =s$ is defined as
\begin{equation}
	R^\perp(s,\hat{n})[F] = \int\limits_ {x \cdot \hat{n} =s} F_{i_1,....i_k}(x) \hat{n}_{i_1} \cdots \hat{n}_{i_k} \mathrm{d}  x
\end{equation} 
For scalars, $k=0$,  this reduces to the standard Radon transform $R$ . 
In the 2D case the normal Radon transform  coincides with the transverse ray transform of Sharafutdinov \cite{sharafutdinov1994integral}, with a slight change of notation, and \cite{kazantsev2004singular} point out that transverse and longitudinal  ray transforms are equivalent in dimension two, with a relabelling of basis vectors. For 3D the normal Radon transform was considered by \cite{polyakova2021singular,svetov2023inversion,svetov2023decomposition}. 
 Considering now the 2D case, at each point $x$, $f$  depends on two unit vectors $\xi_{\mathrm{in}}$ and $\xi_{\mathrm{out }}$ and to apply the normal Radon transform we fix the bistatic angle, that is the angle $\beta$  between  $\xi_{\mathrm{in}}$ and $\xi_{\mathrm{out }}$. The average azimuth direction $\hat{n}$ and the bistatic angle $\beta$ determine the pair $\{ \xi_{\mathrm{in}},\xi_{\mathrm{out }}\}$ with the order ambiguous. We define for a fixed $\beta$, $\tilde{f}(x,\hat{n}) =f(x, \xi_{\mathrm{in}},\xi_{\mathrm{out}})$.  To convert this to a problem in tensor tomography, we choose a fixed natural number $N$, for convenience, even, and approximate $\tilde{f}$ at each point $x$ by a polynomial $\tilde{f}_N$ in the unit vector $\hat{n}$. Consider the terms of degree $k \le N$, for $k$ even, then we can convert this to a homogeneous polynomial of degree $N$ by multiplication by $|\hat{n}|^{(N-k) } =1$.  Similarly for $k$ odd, we can convert to a homogeneous polynomial of odd degree $N-1$. Homogeneous polynomials of degree $N$ can be naturally identified with symmetric tensors of rank $N$. We can uniquely split $R^{\perp}(s,\hat{n})[F]$ into the sum of even and odd parts in $\hat{n}$, and  then attempt to recover rank $N$ even and $N-1$ odd tensor fields separately from that data. 
 
  For the three dimensional case there is an inconvenience that, even if we fix a bistatic angle and the average azimuth unit vector, we need one more scalar  to uniquely  determine the incoming and outgoing direction.  A simple choice is to consider zero bistatic angle, that is the monostatic case $x_T=x_R$, we then just consider a transceiver at a point on a large sphere centred on the scene centre determined by $\hat{n}$. 

In the case were $F = d^N\phi$ for a scalar function $\phi$ and $d$ the symmetric derivative operator, The normal and scalar Radon transform are related by
\begin{equation}
	R^\perp [d^N\phi] = \frac{\partial^N}{\partial s^N}  R[\phi]
\end{equation}
and in this case $\phi$ can be reconstructed by standard inversion methods for the scalar Radon transform. For less smooth $f$ derivatives may have to be understood in the distributional sense. In fact, such a scalar $\phi$ is all that can be reconstructed from normal Radon transform data in dimension 3. Indeed, following the generalized Helmholtz decomposition given in \cite{svetov2023decomposition} the null-space of $R^\perp$ is exactly the tensor fields that are orthogonal to symmetric $N$ th derivatives of scalar fields. This seems rather disappointing from the point of view of anisotropic reconstruction, as we still only recover a scalar at each point with the anisotropy given in terms of its derivatives. For example, a rank 2 tensor field can only be reconstructed as the nearest (in some sense) Hessian matrix of a scalar field.  However, remember that the data we use is from one measurement of the average azimuth in each direction at each point. In practice we might have many different measurement pairs with the same average azimuth.  If we know {\em a priori} information about the directional dependence in $f$ we may be able to use this data to reduce the ambiguity in the reconstruction. 
\begin{figure}[ht]
  \centering%
  \resizebox{0.7\columnwidth}{!}{%
 	\begin{tikzpicture}
		
		\draw[thick, dashed] (0,0) ellipse (5cm and 3cm);
		\coordinate (Rx) at (4,0);   
		\coordinate (Tx) at (-4,0);   
		\coordinate (Target) at (0,3);  

		\node at (Tx) {\includegraphics[width=0.7cm]{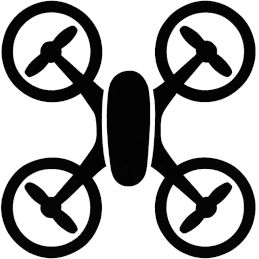}};
		\node[below left=9pt of Tx] {\large $x_T$};	
		\node at (Rx) {\includegraphics[width=0.7cm]{mydrone.png}};
		\node[below right=9pt of Rx] {\large $x_R$};		
	
	    \node at (Target) {\includegraphics[width=1cm]{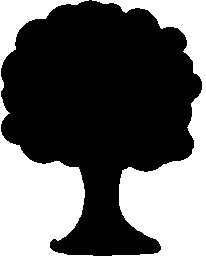}};

		\node[below =0.52cm of Target] {\large Target};
 \pic[
    draw,
    ->,
    thick,
    "$\beta$"{shift={(0,7pt)}},  
    angle eccentricity=1.5,
    angle radius=12mm
] {angle = Tx--Target--Rx};
  \draw[thick, postaction={decorate}, decoration={markings, mark=at position 0.5 with {\arrow[scale=2.0]{>}}}] 
    (Tx) -- (Target);
  \draw[thick, postaction={decorate}, decoration={markings, mark=at position 0.5 with {\arrow[scale=2.0]{>}}}] 
    (Target) -- (Rx);
		
		\node at (0,-3.3) {\large Isochrone (constant time-delay ellipse)};
		
	\end{tikzpicture}
    }
\caption{Isochrones are ellipses in the 2D case with transmitter ($x_T$) and receiver ($x_R$) at the foci. Bistatic angle indicated by $\beta$.\label{fig:isochrone} }
\end{figure}

\section{Preservation of singularities} 
In many applications we are interested in isolated, localized objects. In particular we might have a scene that is largely isotropic with a few small anisotropic objects.  The singular support of a (generalized) function or tensor field is the closure of the set where it fails to be infinitely differentiable. For example, there could be a jump discontinuity where a scene in radar changes type of vegetation or at the boundary of a vehicle or building. The singular support would be precisely that boundary.  The singular support of a delta function is just a point.  For the longitudinal ray transform in dimension 2 we can recover only the {\em solenoidal part}, with the {\em potential part} unknown. For the transverse ray transform we have to relabel basis vectors but the argument is the same. So an important question is to see if we can recover the singular support from just the solenoidal part.  The splitting in to potential and solenoidal parts can be framed in terms of elliptic operators, and these have the property that they preserve singular supports.  In the two dimensional case a tensor $F= G + dV$ where $G$ is solenoidal so $\delta G=0$ where $-\delta$ is the formal adjoint of $d$, so we have 
\begin{equation}\label{eq:deltaF}
    \delta F = \delta dV. 
\end{equation}
Sharafutdinov \cite[lemma 3.3.3]{sharafutdinov1994integral} proves that $\delta d$ is an elliptic partial differential operator.  We can conclude that anything in the singular support of $dV$ is also in the singular support of $F$. This means that while the reconstruction is ambiguous we can expect to locate singularities, and hence discrete objects.  In the three dimensional case where $F= G +d^N \phi$ the adjoint of $d^N$ is $(-1)^N \delta^N$ and we can apply \cite[(6.4.1)]{sharafutdinov1994integral} recursively to show that $\delta^N d^N$
is elliptic and again we can recover the singular support.



\section{Worked example}
Let us consider the rank-2 tensor in dimension two
\[F= \begin{pmatrix} \delta_0 & 0 \\ 0 & -\delta_0 \end{pmatrix}\]
where $\delta_0$ is the delta function at the origin. Then the solenoidal part $G=F-dV$ where 


\begin{equation}
\begin{split}
& \text{\tiny $\displaystyle dV
 = -\frac{1}{2\pi}\operatorname{p.v.}\Bigg ( \frac{1}{r^{2}}
\begin{pmatrix}
3r^4 + \cos(2\theta)(1+3r^4) + \cos(4\theta) & \sin(4\theta)\\
\sin(4\theta) &-3r^4+\cos(2\theta)(1+ 3r^4) - \cos(4\theta)
\end{pmatrix}\Bigg)$}\\
&\hskip3cm \text{\tiny $\displaystyle +\frac{3}{4}
\begin{pmatrix}
-1&0\\[2pt]0&1
\end{pmatrix}
\delta_0$}
\end{split}
\end{equation}


\noindent Here $\mathrm{p.v}$ indicates a distribution  defined by the Cauchy principle value and this formula is determined by solving the equation $\delta d U = \delta F$. We note that the solenoidal part $G$ is no longer localized exclusively at zero but also has other terms that are smooth away from zero and also non-diagonal terms. A  plot of $G$ calculated using this formula are shown in Fig.  \ref{fig:poisson-sol}.

\section{Numerical inversion}

 The solenoidal part $G$ in the above example can also be approximated  by a numerical inversion using a singular value expansion (SVE). In common with X-ray tomography problems \cite{hansen2021computed} reconstruction methods include filtered back projection algorithms, regularized inversion of a linear system and truncated singular value expansion. In the latter case we observe that operators $A$ between one Hilbert space of functions and another that arise as the forward operator in inverse problems often have a families of orthogonal functions  $v_i$ and $u_i$ (singular functions) so that $Av_i = \sigma_i u_i$, for real numbers $\sigma_1 \ge \sigma_2 \ge \cdots>0 $ (singular values),  The $u_i$ can be chosen to span the range of the operator and the $v_i$ the orthogonal complement of the null-space. We can interpret this as $v_i$ being components of the unknown image we wish to reconstruct, but as $\sigma_i$ decreases they make smaller contributions to the measured data. The more rapidly the singular values decrease, the more ill-posed the problem. To solve $A f= g$ approximately we choose an M so that $\sigma_i$ is around the level of our measurement and modelling error for $i>M$, and approximate $f$ by its truncated singular value expansion
\[
f\approx \sum\limits_{i=1}^M \frac{\langle g, u_i \rangle v_i}{\sigma_i}.
\]
While not many operators in imaging inverse problems have explicitly known SVEs we are fortunate that the normal Radon transform for symmetric tensor fields on the unit disk does \cite{kazantsev2004singular}. The tensor fields corresponding to the $v_i$ have components that are Zernike disk functions, while in the correct angular coordinates the data is expanded in combinations of Fourier basis functions. For the 3D case \cite{polyakova2021singular} gives the SVE for the normal Radon transform of symmetric tensor fields on the unit ball.



\begin{figure}[htbp]
  \centering
  \hspace*{-0.5cm}
  \includegraphics[width=.5\textwidth]{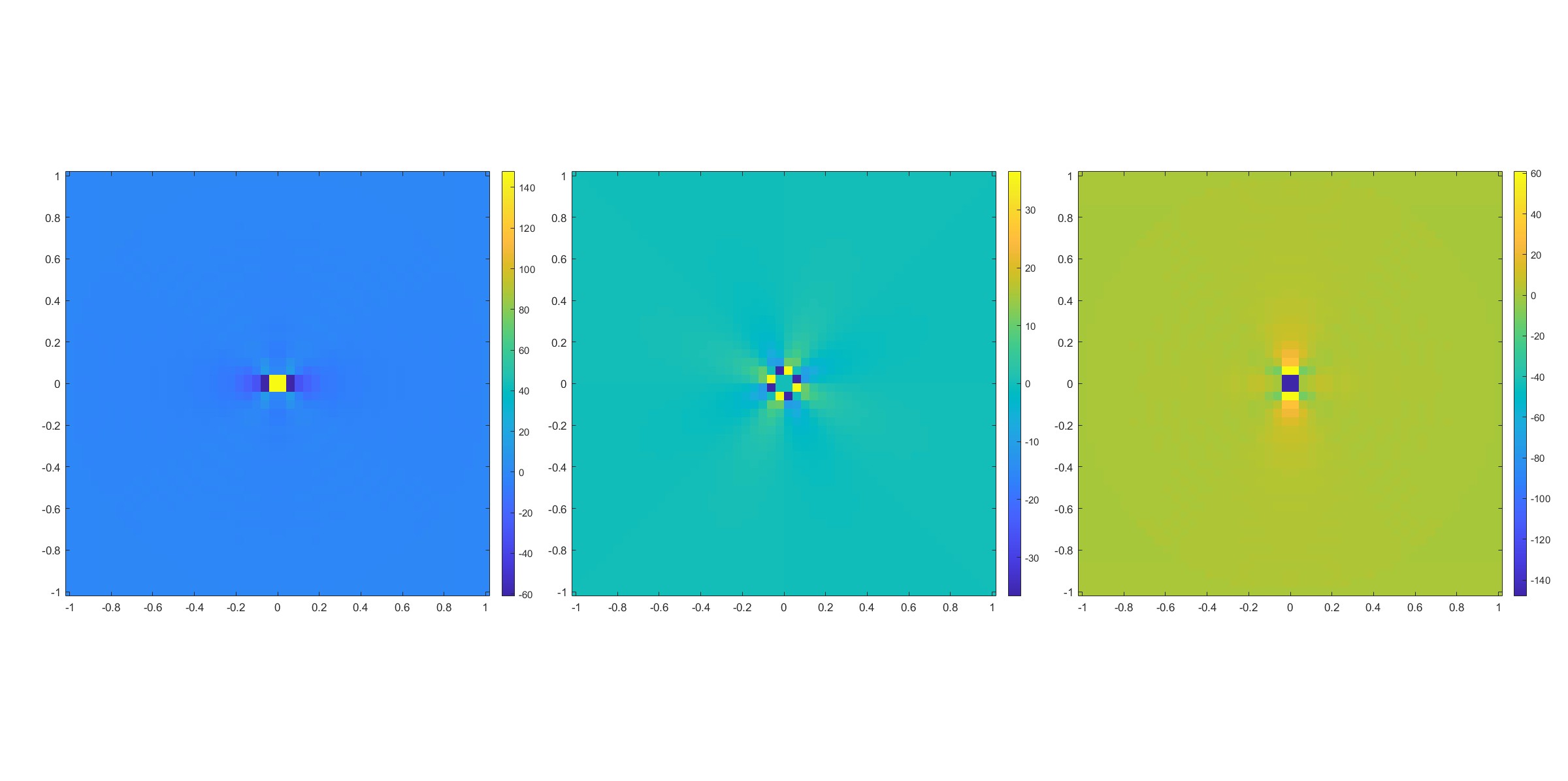}
  \captionsetup{skip=-20pt} 
  \caption{Solenoidal components $G$ from solving $\delta d U = \delta F$ with $G_{11}$ on the left, $G_{12}$ in the centre and with $G_{22}$ on the right. }
  \label{fig:poisson-sol}
\end{figure}



\begin{figure}[htbp]
  \centering
  \hspace*{-0.5cm}
  \includegraphics[width=.5\textwidth]{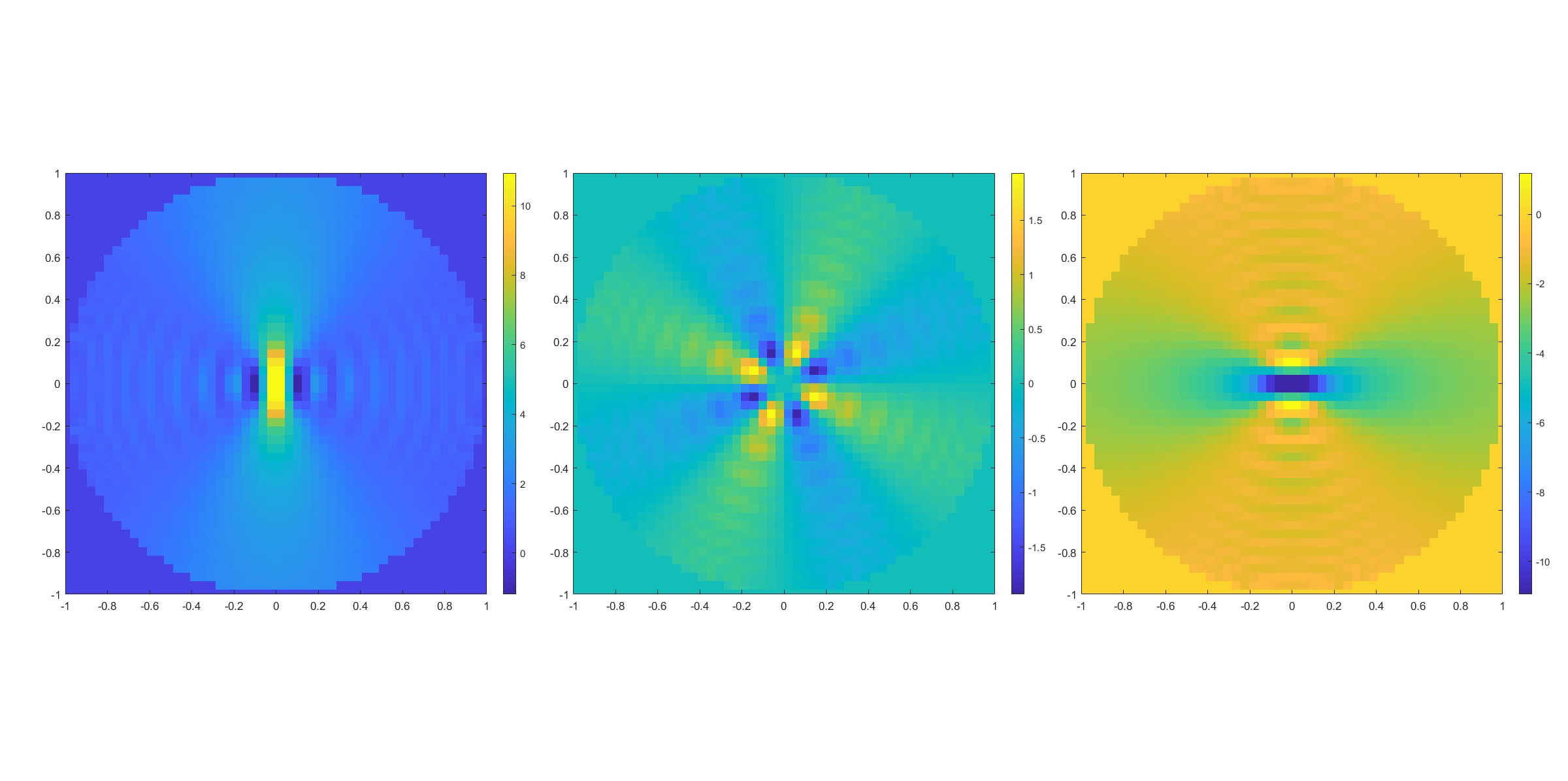}
  \captionsetup{skip=-20pt} 
  \caption{Reconstructed solenoidal components via the Truncated SVE of an isolated object for \(N=50\) and \(K=40\), with $a_0$ on the
left, $a_1$ in the middle and with $a_2$ on the right}
  \label{fig:KB-sol}
\end{figure}

Now we will compare Fig.  \ref{fig:poisson-sol} against determination of $G$ using the truncated SVE described in \cite{kazantsev2004singular}, in which the singular functions are parameterized by two integers $K$ the angular frequency and $N$ the order of the radial polynomial. We implement the truncated SVE for a mesh grid $N=50$ and $K=40$ as the truncation parameter where $K<N$ with its components presented in Fig.  \ref{fig:KB-sol}.



\section{Conclusions and further work}

We observe in Fig.  \ref{fig:KB-sol} that, under truncated SVE, `ringing' oscillations appear in the reconstructions;   this is consistent with the Gibbs phenomenon that occurs when representing functions with singularities by finite Fourier series 

In the case of the airborne far-field setting in which the transmitter and receivers are confined to the sky (upper half–space) and only a restricted set of line (or plane) normals is sampled. The set of admissible look–directions is therefore a hemisphere, so for multistatic measurements the accessible pairs of directions form a subset of the product \(\mathbb S^{2}_{+}\times\mathbb S^{2}_{+}\) (the Cartesian product of a hemisphere with itself).


A practical use case is the detection of isolated objects against an isotropic background, though we are specifically concerned with detecting isolated anisotropic objects  with a localized singular support. This is illustrated by our worked example for the deviatoric delta case for a rank 2-tensor field. 


\begin{acks}
The authors would like to thank the Isaac Newton Institute for Mathematical Sciences, Cambridge, for support and hospitality during the programme RNT where work on this paper was undertaken. This work was supported by EPSRC grant nos EP/Z000580/1, EP/V007742/1, and an iCASE award with DSTL. WL completed part of this work while a Visiting Fellow at Clare Hall Cambridge.
\end{acks}

\balance

\nocite{*}
\bibliography{anisott}
\bibliographystyle{iet}


\end{document}